\documentclass[11pt]{article}
\usepackage{amsfonts,pdfsync}
\usepackage{amsthm}

\usepackage{amsmath}
\usepackage{graphicx}
\usepackage{comment}
\usepackage[usenames,dvipsnames]{xcolor}
\usepackage{enumerate}
\usepackage{extarrows}
\usepackage{hyperref}
\usepackage[round, authoryear]{natbib}
\usepackage{stackrel}
\usepackage{centernot}
\usepackage{caption}
\usepackage{subcaption}
\usepackage{fancyhdr}

\let\quoteOLD\quote
\def\quote{\quoteOLD\small}

\definecolor{labelkey}{cmyk}{0,0.8,1,0.5}
\definecolor{refkey}{cmyk}{0,0.8,1,0.5}

\newtheorem{theorem}{Theorem}[section]
\newtheorem{example0}{\sc Example}[subsection]

\newtheorem{proposition}{Proposition}

\newtheorem{remark}{Remark}

\numberwithin{equation}{section}
\numberwithin{theorem}{section}
\numberwithin{corollary}{section}
\numberwithin{proposition}{section}
\numberwithin{lemma}{section}
\numberwithin{definition}{section}
\numberwithin{remark}{section}

\makeatletter
\def\th@newremark{\th@remark\thm@headfont{\bfseries}}
\makeatletter

\def\boxit#1{\vbox{\hrule\hbox{\vrule\kern6pt
          \vbox{\kern6pt#1\kern6pt}\kern6pt\vrule}\hrule}}

\newcommand{\sid}[1]{{\color{black} #1}}
\newcommand{\res}[1]{{\color{black} #1}}

\newcommand{\pibar}{\overline{\Pi}}
\newcommand{\pibarinv}{\overline{\Pi}^{\leftarrow}}

\newcommand{\cF}{{\cal F}}

\newcommand{\Fbar}{\overline{F}}

\newcommand{\Iinv}{I^\leftarrow}

\newcommand{\Gbar}{\overline{G}}

\newcommand{\veps}{\varepsilon}

\newcommand{\topr}{\stackrel{\mathrm{P}}{\longrightarrow}}
\newcommand{\todr}{\stackrel{\mathrm{D}}{\longrightarrow}}

\newcommand{\eqdr}{\stackrel{\mathrm{D}}{=}}

\newcommand{\eqd}{\stackrel{{\mathrm D}}=}

\newcommand{\R}{\Bbb{R}}

\newcommand{\N}{\Bbb{N}}

\newcommand{\rmd}{{\rm d}}

\newcommand{\halmos}{\quad\hfill\mbox{$\Box$}}

\newcommand{\BB}{\mathbb{B}}
\newcommand{\DD}{\mathbb{D}}

\newcommand{\NN}{\mathbb{N}}

\newcommand{\XX}{\mathbb{D}}

\newcommand{\AAA}{{\cal A}}

\newcommand{\wt}{\widetilde}

\newcommand{\dto}{\downarrow}

\newcommand{\EEEE}{\mathfrak{E}}

\newcommand{\be}{\begin{equation}}
\newcommand{\ee}{\end{equation}}
\newcommand{\bea}{\begin{eqnarray}}
\newcommand{\eea}{\end{eqnarray}}
\newcommand{\bean}{\begin{eqnarray*}}
\newcommand{\eean}{\end{eqnarray*}}
\newcommand{\ben}{\begin{equation*}}
\newcommand{\een}{\end{equation*}}
\newcommand{\ba}{\begin{aligned}}
\newcommand{\ea}{\end{aligned}}

\def\nexto{\kern -0.54em}

\newcommand{\PP}{\textbf{\rm P}}
\newcommand{\EE}{\textbf{\rm E}}

\newcommand{\BN}{\mathcal{BN}}

\begin{document}

\title{\bf Ratios of Ordered Points of Point Processes with Regularly Varying Intensity Measures}
\author{Yuguang Ipsen, Ross Maller and Sidney Resnick\thanks{Ross Maller was
partially supported by ARC Grants DP1092502  and DP160103037. 
Sidney Resnick was partially supported by US ARO MURI grant W911NF-12-1-0385 and ARC Grant DP160103037. Hospitality and space from the Australian National University in February 2017 is \res{acknowledged.} \newline
Email: yuguang.ipsen@anu.edu.au; ross.maller@anu.edu.au; sir1@cornell.edu}}


\maketitle

\begin{abstract}
We study limiting properties of ratios of ordered points of point processes whose intensity measures have regularly varying tails,  giving a systematic treatment which  points the way to 
 ``large-trimming" properties of extremal processes and a variety of applications.
Our point process approach facilitates a connection with the negative binomial process of \cite{Gregoire1984}
 and consequently to  certain generalised versions of the Poisson-Dirichlet distribution.
\end{abstract}

\section{Introduction}\label{intro}

Recent work on ratios of \sid{ordered Poisson points}  and
\sid{ordered} jumps of
stable subordinators  and other L\'evy processes due to \cite{KeveiMason2014} and the present
authors in \cite*{bfm2016},  
\cite*{stableJump1} and \cite*{BMR2016}
placed  an emphasis  on limiting properties of those  ratios, and on ``trimmed" versions of the \sid{ process
  generating the points}, which may have been a subordinator or a more general L\'evy process. 
  

Our aim in \sid{this} paper is to give a systematic treatment of 
the limiting behaviour of ratios of ordered Poisson points. As is natural, we take a point process approach and make special
  connection with the negative binomial process whose relevance in the
  present context was brought out in \cite{trimfrenzy}.  
This connection via ratios of points enabled the construction of a generalised kind of
Poisson-Dirichlet distribution which can be added to the repertoire of
available models for data analytic purposes.  

\sid{A related topic is the behaviour of two dimensional Poisson points
  ordered by the second component
  when the $r$ highest points are deleted.} 
\sid{Such processes} were explored in \cite*{BMR2016},
and the present results  \sid{provide}  impetus for further investigations of this kind. 

The paper is structured as follows.
In Section \ref{s2} we set up  the  point processes to be studied in 
 Section \ref{s3}, notably a Poisson point process $\DD$ on $\R^+ = (0, \infty)$, 
 and subsidiary  point processes  $\XX^{(n)}_t$ and 
  $\XX_t^{(r,r+n)}$ 
 consisting of ratios of the ordered points in $\DD$, where the ordering is by magnitude up till a given time $t>0$, and the normalisation is by the $n$th largest point. 

The tail of the canonical measure for the points is assumed to be regularly varying of index $-\alpha$, $\alpha>0$, at 0.
Under this assumption, Theorem \ref{preT_2} in Section \ref{s3} proves the weak convergence of $\XX_t^{(r,r+n)}$, as $t \dto 0$, 
to a limit comprised of a sum of independent point processes on $(0,\infty)$. The first component of the sum 
represents the joint limiting distribution of
 ratios larger than 1 of points in $\DD$, conveniently expressed as the distribution of the order statistics of certain i.i.d. (independent, identically distributed) random variables (rvs);
  and the second component is a negative binomial point process, representing the  limiting distribution of the (infinitely many)  ratios smaller than 1.

Further, in Section \ref{s4}, we mention some interesting corollaries of  Theorem \ref{preT_2}, 
stated as separate propositions, 
and in Section \ref{s5} prove a converse result (Theorem \ref{conv}) 
to the effect  that convergence in distribution of ratios (larger or smaller than 1) implies regular variation of
 the tail of the canonical measure for points in $\DD$.
We conclude in Section \ref{s6} with some history relating to antecedents of these results in the  literature of order statistics of i.i.d.
rvs, which can be used to suggest further explorations in that area.

\section{Poisson Point Processes and Ratios of Ordered Points}\label{s2}
In this section we set up the  point process framework we will use.
\res{Let $\N=\{1,2,\ldots\}$ and $\N_0=\{0,1,2,\ldots\}$.
Suppose $\Pi$ is a Borel measure on $(0,\infty)$, locally finite at
infinity. The measure $\Pi$ has finite-valued tail function $\pibar:(0,\infty)\to (0,\infty)$, defined by 
\ben  
\pibar(x):= \Pi\{(x,\infty)\}, \ x>0,
\een
a right-continuous, non-increasing function.
Assume throughout that $\Pi\{(0,\infty)\}=\pibar(0+)=\infty$, so there are infinitely many non-zero points of
$\DD$ in any right neighborhood of $0$.
 Let
\ben  
\pibar^\leftarrow(x)=\inf\{y>0: \pibar (y) \le x\},\ x>0,
\een
be the right-continuous inverse of $\pibar$. }
\res{With $\delta_{x}$ denoting a point mass at $x\in(0,\infty)$, let
\be\label{Ddef}
\DD_t=\sum_{j>0}\delta_{\Delta_t^{(j)}},
\ee
a Poisson point process on $(0,\infty)$ with intensity measure
$t\Pi(\rmd x)$, 
where the points are written in decreasing order, possibly with ties,}
$$\infty>\Delta^{(1)}_t\ge {\cdots}\ge \Delta^{(r)}_t\geq \dots >0 .$$

A representation detailed in \cite*{bfm2016} 
shows how to construct all $\DD_t$  processes on the same space.
Since $\pibar(0+)=\infty$, all $\Delta^{(r)}_t$ are positive and $\Delta^{(r)}_t\dto 0$ a.s. as $t\dto 0$ for $r\in\N$.
Let $(\EEEE_i)$ be an i.i.d. sequence of exponentially distributed random variables with common parameter $\EE\EEEE_i=1$.
Then $\Gamma_r:=\sum_{i=1}^r\EEEE_i$ is a Gamma$(r,1)$ random
variable, $r\in\N$,
and $\{\Gamma_r, r\ge 1\}$ are the points of a homogeneous, unit rate Poisson
process on $\R^+$.
\sid{The representation is}
\be\label{ch}
\big\{\Delta_t^{(i)}\big\}_{ i \ge  1}
\eqdr \big\{\pibarinv(\Gamma_i/t)\big\}_{ i \ge 1},\ t>0.
\ee
For earlier and  related representations  consult
   \cite{LeP1980a, LeP1981}; \cite*{LePZinn1981};
   \cite{samtaqqu}, \res{p. 21, 30;     \cite{res87}, Ex. 3.38, p.139;
   \cite{resnick:1986b}, Sect. 2.4; and \cite{ferguson:klass:1972}.}

Write
\ben
\PP(\Gamma_r\in\rmd x)= \frac{x^{r-1} e^{-x}\rmd x}{\Gamma(r)} {\bf 1}_{\{x>0\}},\ r\in\N,
\een
for the density of $\Gamma_r$, 
which should not be confused with the Gamma function,
$\Gamma(r)=\int_0^\infty x^{r-1}e^{-x}\rmd x$, $r>0$.
A beta random variable ${\rm B}_{a, b}$ on $(0,1)$  with parameters $a, b>0$ has density function
\[ f_B(x) = \frac{\Gamma(a+b)}{\Gamma(a)\Gamma(b)} x^{a-1}(1-x)^{b-1} = \frac{1}{B(a,b)}x^{a-1}(1-x)^{b-1},\  0<x<1.
\] 
Thus
\be\label{incombet}
\PP\left({\rm B}_{a, b}\le x\right) =  \frac{1}{B(a,b)}\int_0^x y^{a-1}(1-y)^{b-1} \rmd y
:= B(a,b;x),\ 0<x<1,
\ee
where $ B(a,b;x)$ is the incomplete Beta function.
%

We are interested in the convergence behaviour of ratios of the order statistics $\Delta^{(r)}_t$, as $t\dto 0$.
The basic assumption is the regular variation of the tail function $\pibar(x)$.
Write $RV_0(\beta)$ (resp. $RV_\infty(\beta)$) for the real-valued functions  regularly varying at 0 (resp, infinity) with index $\beta$. 
We have  $\pibar(x)\in RV_{0/\infty}(-\alpha)$, $0\le \alpha\le \infty$, iff
\ben
\lim_{x\to 0/\infty}\frac{\pibar(x\lambda)}{\pibar(x)}=\lambda^{-\alpha},\
{\rm for}\ \lambda>0.
\een
Interpret
$\lambda^{-\infty}=0.{\bf 1}_{\{\lambda>1\}}+1.{\bf 1}_{\{\lambda=1\}}
+\infty.{\bf 1}_{\{\lambda<1\}}$ and 
$1/0\equiv \infty$. 
From \citet*[p.28-29]{BGT87} we know that $\pibar(x)\in RV_{0/\infty}(-\alpha)$ iff
  $\pibarinv(x)\in RV_{\infty/0}(-1/\alpha)$.
\sid{The slowly varying functions at $0$ or $\infty$ are denoted  $RV_{0/\infty}(0)$}  and $RV_{0/\infty}(\infty)$ are the rapidly varying functions at $0$ or $\infty$.

When $\pibar(\cdot) \in RV_0(-\alpha)$ with $0\le \alpha\le \infty$ or, equivalently,  $\pibarinv(\cdot) \in RV_\infty(-1/\alpha)$,  we have the easily verified convergence
(with the interpretation as above when $\alpha=0$ or $\alpha=\infty$)
\be\label{1b}
t\pibar(u\pibarinv(y/t))
\sim
\frac{\pibar(uy^{-1/\alpha}\pibarinv(1/t))} 
{\pibar(\pibarinv(1/t))}
 \to u^{-\alpha} y \ \text{as} \ t \dto 0,\ {\rm for\ all}\ u, y > 0.
\ee

\section{Ratios of Ordered Points}\label{s3}
In this section we give a general result
for the point processes of ratios of ordered points of $\DD_t$.
Fix $r \in \N_0$, $n \in \N$ and $t>0$.
Define the point processes  on $(0,\infty)$:
 \be\label{PP0}
\XX^{(n)}_t := \sum_{i \ge n+1} 
\delta_{\{ \Delta^{(i)}_t/\Delta^{(n)}_t \,\} }
\ee 
and 
\be\label{pp}
\XX^{(r, r+  n)}_t :=  \sum_{i \ge r+1}  \delta_{\{\Delta^{(i)}_t/\Delta^{(r+n)}_t \}}
 = \sum_{i = 1}^{n-1} \delta_{\{ \Delta_t^{(r+i)}/\Delta_t^{(r+n)}\}}+1 +\XX_t^{(r+n)}.
 \ee
Conditionally on $\{ \Delta_t^{(n)}=z\}$, $z > 0$, the points
$(\Delta_t^{(i)})_{i\ge n+1}$ comprise a Poisson point process with
intensity measure $\sid{t}\Pi$ restricted to $(0,z)$. 
Thus, the Laplace functional of $\XX_t^{(n)}$, conditional on $\{ \Delta_t^{(n)}=z\}$, is
\begin{align}\label{Lapf}
\EE \big(e^{-\XX_t^{(n)}(f)} \, \big| \,  \Delta_t^{(n)}=z \big) 
=&\, \EE\Big(\exp\Big(-\int_{0<x<1} f(x) \XX_t^{(n)}(\rmd x)\Big)
\Big| \, \Delta_t^{(n)} = z\Big)\nonumber \\
=& \, \exp\Big(-t\int_{0<x<1}(1-e^{-f(x)})\Pi(z\rmd x)\Big),
\end{align}
where $f\in \cF_+$, the nonnegative measurable functions on $\R^+$.

Let $\mathbb{S}$ be a Poisson point process \sid{on $(0,\infty)$} with
  intensity measure $\Lambda(\rmd x) = \alpha x^{-\alpha-1}{\bf
    1}_{\{x>0\}}$,  
$ \alpha>0$, represented as $\mathbb{S} = \sid{\sum_{i\ge 1}
  \delta_{\Gamma_i^{-1/\alpha}}}$. When $0<\alpha<2$ we can interpret
\sid{$\Gamma_i^{-1/\alpha}$} as the $i$th largest jump of a stable process
$(S_t)_{0<t\le 1}$ with L\'evy measure $ \Lambda$, but we allow any
$\alpha>0$. 
Analogous to \eqref{PP0}, define
\be\label{defBn}
\BB^{(n)} = \sum_{i\ge n+1} \sid{\delta_{\{(\Gamma_i/\Gamma_n)^{-1/\alpha}\}},}\ n\in\N.
\ee
   
The  point process in \eqref{defBn} has Laplace functional at $f$ equal to 
\be\label{NBlap}
\EE(e^{-\BB^{(n)}(f)}) = \Big(1 + \int_{0}^1 \big(1 - e^{-f(x)} \big) \Lambda(\rmd x) \Big)^{-n}.
\ee
$\BB^{(n)}$ is the negative binomial point process with base measure $\Lambda^*(\rmd x) = \Lambda(\rmd x){\bf 1}_{0<x<1}$,  denoted by $\BN(n,  \Lambda^*)$, in the notation of \cite{Gregoire1984}.

The next theorem shows the weak convergence (denoted by `$\todr$') of $\XX_t^{(r,r+n)}$ as $t \dto 0$ to a limit comprised of independent components of  $\BB^{(r+n)}$ and a mixture of beta random variables. 

\begin{theorem}\label{preT_2} 
Suppose $\pibar(\cdot) \in RV_0(-\alpha)$, $0<\alpha<\infty$, and $f\in \cF_+$.
 Fix $n,\sid{r}  \in \N$.
Then 
\begin{enumerate}[{\rm (i)}]
\item \sid{In the space of point measures $M_p(0,\infty)$ with the
    vague topology, as $t\dto 0$,
\begin{equation}\label{Mt0}
\XX^{(r,r+n)}_t
\todr  \sum_{j=1}^\infty \delta_{\{(\Gamma_{r+j}/\Gamma_{r+n})^{-1/\alpha}\}} =
\sum_{j=1}^{n-1}
\delta_{\{(\Gamma_{r+j}/\Gamma_{r+n})^{-1/\alpha}\}} 
+\delta_{1} +
\BB^{(r+n)}.
\end{equation}
The limit has Laplace functional at $f$ equal to }
\begin{equation}\label{e:onemore}
\EE\big( e^{-f(J(B_{r,n}^{1/\alpha}))} \big)^{n-1}\,
 e^{-f(1)}\,
\EE\big(e^{-\BB^{(r+n)}(f)} \big),
\end{equation}
where, for each $u\in (0,1)$, $J(u)$ has distribution
\be\label{Jdis0}
\PP(J(u)\in\rmd x)= \frac{\Lambda(\rmd x) 
{\bf 1}_{\{1<x<1/u\}}}{1-u^{\alpha}},\ x>0,
\ee
$B_{r,n}$ is a Beta$(r,n)$  random variable independent of $J(u)$,
\res{and the third factor on the right of \eqref{e:onemore} is
    determined from  \eqref{NBlap}.  }

\item For $r=0$,
\be\label{Mt1}
\lim_{t\dto 0}
\EE\big( e^{- \XX^{(0,n)}_t (f)}\big)\,
= \EE\big( e^{-f(L)} \big)^{n-1} e^{-f(1)}\,
 \EE\big(e^{-\BB^{(n)}(f)} \big),
\ee
where 
$L$ is a  random variable  with distribution
\be\label{Ldis}
\PP(L \in \rmd x)= \Lambda(\rmd x) {\bf 1}_{\{x > 1\}}.
\ee
\end{enumerate}
\end{theorem}

\medskip\noindent{\bf Proof of  Theorem \ref{preT_2}:}\
(i)\ Using the representation in \eqref{ch}
and the fact that $\pibarinv \in RV_{\infty}(-1/\alpha)$, we immediately get, as $t\dto 0$, with almost sure convergence,
\begin{align}\label{rat2}
  \Bigg( \frac{\Delta_t^{(r+j)} }{\Delta_t^{(r+n)} };\, j \geq 1 \Bigg)
  \eqd \Bigg( \frac{\pibarinv(\Gamma_{r+j}/t) }{\pibarinv(\Gamma_{r+n}/t)};\, j \geq 1 \Bigg) 
    \to  
  \Bigg( \Bigl(\frac{\Gamma_{r+j} }{\Gamma_{r+n} }
   \Bigr)^{-1/\alpha}; \,j\geq 1\Bigg)
\end{align}
for each $r,n\in\NN$. 
By separating the ratios in the limit process into those bigger than 1, equal to 1, or smaller than 1, we get the form in \eqref{Mt0}. 

The points in the  two limit point processes in \eqref{Mt0} occur in non-overlapping regions, so,
conditionally on $\Gamma_{r+n}$,  they are independent of each other. In fact, by the algebraic properties of gamma distributions, i.e., $\Gamma_r/\Gamma_{r+n} \eqd B_{r,n}$ with $B_{r,n}$ independent of $\Gamma_{r+n}$, these components are also unconditionally independent. Thus the Laplace transform can be given  in the product form of \eqref{e:onemore}. Next we will derive the Laplace functional for each component separately.

For ratios bigger than $1$, we  note from properties of a  homogeneous
Poisson process that, \res{conditionally on $\Gamma_r/\Gamma_{r+n} =
B_{r, n}=s$,  
$$\frac{\Gamma_{r+1}}{\Gamma_{r+n}} ,\dots, \frac{\Gamma_{r+n-1}}{\Gamma_{r+n}}
$$
are the order statistics of a uniform sample of size $n-1$ on $(s,1)$
and the unordered sample has representation 
$$s+(1-s)U_j; \; j=1\dots,n-1,$$
where $U_1,\dots,U_{n-1}$ are iid uniform on $(0,1)$. Thus for 
$f \in \cF_+$ , 
\begin{align}\label{313}
\EE \exp \Big( {-\sum_{j = 1}^{n-1}
  f\Big(\Big(\frac{\Gamma_{r+j}}{\Gamma_{r+n}} \Big)^{-1/\alpha}}
  \Big) & |B_{r,n}=s\Big) 
= \EE \exp \big\{
 {-\sum_{j = 1}^{n-1} f \bigl((1-s)U_j+s)^{-1/\alpha} \bigr) } \bigr\}  \nonumber \\
=&\Bigl(
\int_0^1  \exp \big\{
 {- f \bigl((1-s)u+s)^{-1/\alpha} \bigr) } \bigr\} du \Bigr)^{n-1}, \nonumber \\
\intertext{and setting  $y=((1-s)u+s)^{-1/\alpha}$ gives}
=&\Bigl(\int_1^{s^{-1/\alpha}}  e^{-f(y)} \alpha y^{-\alpha -1}
   \frac{dy}{1-s}\Bigr)^{n-1}.
\end{align} 
Take expectations in \eqref{313} to get
$$
\EE \exp \Big( {-\sum_{j = 1}^{n-1}
  f\Big(\Big(\frac{\Gamma_{r+j}}{\Gamma_{r+n}} \Big)^{-1/\alpha}}
  \Big) \Big) 
=\EE  \Bigl(\int_1^{B_{r,n}^{-1/\alpha}}  e^{-f(y)} \alpha y^{-\alpha -1}
   \frac{dy}{1-B_{r,n}}\Bigr)^{n-1}
  $$
which gives \eqref{Jdis0}.}

Next we compute the intensity measure of the limit point process with ratios less than 1, that is, the process $\DD_t^{(r+n)}=\sum_{j\ge n+1} \delta_{\{(\Gamma_{r+j}/\Gamma_{r+n})^{-1/\alpha}\}} $ in \eqref{pp}. 
Conditionally on $\Gamma_{r+n}$, the process
$\sum_{j\ge n + 1} \delta_{\{\Gamma_{r+j}/\Gamma_{r+n}\}}$ 
is a Poisson process with mean measure $\Gamma_{r+n} \rmd x$, where $\rmd x$ is the Lebesque measure. Then the image measure of $\Gamma_{r+n} \rmd x$ under the map $T:x\mapsto  x^{-1/\alpha}$ is $\Gamma_{r+n} \Lambda(\rmd x) $. Hence for any nonnegative measurable function $f$, 
\begin{align}\label{e:NBlap}
& \EE \exp \Big( -\sum_{j\ge n+1} f\Big(\Big(\frac{\Gamma_{r+j}}{\Gamma_{r+n}} \Big)^{-1/\alpha} \Big) \Big) \nonumber \\
 &= \int_{y>0} \exp\Big(-\int_{0<x <1}\big(1-e^{-f(x)}\big) y \Lambda(\rmd x) \Big)  \PP(\Gamma_{r+n} \in  \, \rmd y) \nonumber \\
 &=  \Big( 1 + \int_{0}^1 \big(1-e^{-f(x)}\big) \Lambda(\rmd x) \Big)^{-r-n}.
\end{align}
Referring to \eqref{NBlap}, this is the Laplace transform of a negative binomial point process $\BN(r+n, \Lambda(\rmd x){\bf 1}_{0<x<1})$ at $f$.

(ii)\ ($r=0$) \ The proof of  \eqref{Mt1} is very similar. The
treatment for ratios smaller than or equal to 1 is exactly the
same. 
\halmos
\begin{remark}\label{rem1}
{\rm 
The first component on the RHS of the limit in \eqref{e:onemore} 
shows that, after deleting the $r$
largest points,  the sum  \sid{
$\sum_{i=r+1}^{r+n-1}(\Gamma_i/\Gamma_{r+n} )^{-1/\alpha}$}
has the distribution of a sum of i.i.d.  random variables, once we condition on a $B_{r,n}^{1/\alpha}$  random variable. 
The third component on the RHS of  \eqref{e:onemore}
 is the negative binomial point process $\BB^{(r+n)}$ with base measure $\Lambda$.
So we have the nice representation resulting from the decomposition of the original process  into parts including ratios smaller than 1 and greater than 1.

Ratios of jumps of stable subordinators also featured prominently in the work of 
\cite{PY1997}.
Much subsequent related research involved Poisson-Dirichlet distributions and their involvement in fragmentation and coalescence problems; see  \cite{Bertoin2006} and references therein.
An early influential paper was \cite{Kingman1975}.
The resulting processes have found wide application in a variety of applied areas ranging from Bayesian statistics to models for species diversity;  see for example the list in \citet[Sect.1]{PY1996}.

When $0<\alpha<1$, $r = 0$, similar results were obtained  in Lemma 24 of \cite{PY1997} but without explicit reference being made to the negative binomial point process of \cite{Gregoire1984}.
 Our result allows the bigger range of $\alpha$, $\alpha>0$,  and generalises to point processes with  intensity measures whose tails are regularly varying, rather than dealing only with jumps of subordinators. 
 In general, in our scenario, the points of the limiting process may not be summable. As a special case, for example, we deal elsewhere 
(in \cite{stableJump1}) 
 with L\'evy processes in the domain of attraction of a stable process with index $\alpha\in(1,2)$; compensating the process is then essential.
 
}
\end{remark}
In the next section we draw out some ramifications of Theorem \ref{preT_2}.

\section{Corollaries, Special Cases and Further Results}\label{s4}
Theorem \ref{preT_2} is \res{expressed as convergence of point
  processes.}
\res{In this section we express the theorem in a different form in order to 
facilitate comparisons with earlier results in the literature; we also
extend the result 
to the  $\alpha=0$ or $\alpha=+\infty$ cases and
consider limits of conditional distributions.}

The discussion  is again conveniently divided into parts covering ratios smaller than  or greater than 1. 
Define the ratio
\be\label{Cw}
W_{r,n}(t):=
\frac{\Delta^{(r+n)}_t}{\Delta^{(r)}_t},\ r,n\in\N,\ t>0.
\ee

\begin{proposition}[Ratios smaller than 1]\label{cor2a}
Suppose $\pibar(\cdot) \in RV_0(-\alpha)$, 
$0\le \alpha\le \infty$. 

{\rm (i)}\ Suppose $0<\alpha<\infty$. Then, for each $n \in \N$,
\be\label{weak}
\XX_t^{(n)} \todr  \BB^{(n)},\ {\rm  as}\ t \dto 0,
\ee
where 
 $\BB^{(n)}$ is distributed as $\BN(n, \Lambda)$ with Laplace functional as in \eqref{NBlap}. 
 
{\rm (ii)}\
 As $t \dto 0$, for each $r,n \in \NN$,
\be\label{asymIndep}
\Bigg(\frac{\Delta^{(r+1)}_t}{\Delta^{(r)}_t}, \ldots,  \frac{\Delta^{(r+n)}_t}{\Delta^{(r+n-1)}_t}\Bigg) 
\todr \big(Y_{r}, \ldots, Y_{r+n-1}\big),
\ee
where $Y_{k}$, $k \in\N$, are mutually independent random variables with $Beta(k\alpha, 1)$ distributions. When $\alpha=0$ or $\alpha=\infty$, \eqref{asymIndep} remains true with each $ Y_{k}$ equal to 0 or  with each $ Y_{k}$ equal to 1, respectively.

{\rm (iii)}\ When  $0<\alpha<\infty$, 
$W_{r,n}(t)$ in \eqref{Cw} has limiting distribution as $t\dto 0$ that of 
\be\label{Wg}
 \prod_{i = 1}^n Y_{r+i-1}\eqdr
\left(\frac{\Gamma_r}{\Gamma_{r+n}}\right)^{1/\alpha}
\eqdr B_{r,n}^{1/\alpha} 
 =: W_{r,n},
\ee
where $W_{r,n}$ has density
\be\label{prodBeta2}
f_{W_{r,n}}(w) = \frac{(1-w^\alpha)^{n-1}\alpha w^{\alpha r- 1}}{B(r,n)},\ 0 < w < 1.
\ee
\end{proposition}

\medskip\noindent{\bf Proof of Proposition \ref{cor2a}:}\
(i)\ The convergence in \eqref{weak} is an immediate consequence of \eqref{Mt1}.

(ii)\ 
When  $0<\alpha<\infty$ the convergence in \eqref{asymIndep} follows immediately from \eqref{rat2}. When $\alpha=0$ or $\alpha=\infty$, \eqref{rat2} remains true with the appropriate interpretations as outlined in the discussion leading to \eqref{1b}.

(iii)\  \res{Equation} \eqref{Wg} \res{is implied by} \eqref{asymIndep} and the density in \eqref{prodBeta2} is easily calculated. 

\halmos

\begin{remark}\label{remKV}
{\rm 
Treated as ratios of ordered  jumps of a subordinator, \cite{KeveiMason2014} proved the case $n=1$ in \eqref{asymIndep}, among other results 
comparing the magnitudes of ordered  jumps of a subordinator with the
 magnitude of the subordinator itself.
  Proposition \ref{cor2a} is a multidimensional version of their Theorem 1.2, with the $(\Delta_t)$ treated as points  in $\DD$, in their own right.
  (They also proved converse results; see Section \ref{s5}.)
 }
\end{remark}

\begin{proposition}[Ratios greater than 1]\label{cor1}
Suppose $\pibar(\cdot) \in RV_0(-\alpha)$ with $0\le \alpha\le\infty$. 
Take $x_k\ge 1$ for  $0\le k\le n-1$, $n=2,3,\ldots$, $r\in\N$ and $z>0$. 

{\rm (a)}\  Assume $0<\alpha <\infty$ and $z > 0$. 

{\rm (i)}\ Then, for  $0<u<1$,
\bea\label{13b}
&&
\lim_{t\dto 0}
\PP\bigg(\frac{\Delta_t^{(r+k)} }{\Delta_t^{(r+n)}}>x_k,\ 0\le k\le n-1\bigg|W_{r,n}(t)=u, \, \Delta_t^{(r+n)}= \pibarinv(z/t)\bigg)\cr
&&\cr
&&=
\lim_{t\dto 0}
\PP\bigg(\frac{\Delta_t^{(r+k)} }{\Delta_t^{(r+n)}}>x_k,\ 0\le k\le n-1\bigg|W_{r,n}(t)=u\bigg)\cr
&&\cr
&&=
{\bf 1}_{\{u<x_0^{-1}\}}
\PP\big(J_{n-1}^{(k)}(u)>x_k,\ 1\le k\le n-1\big),
\eea
where
$J_{n-1}^{(1)}(u)\ge J_{n-1}^{(2)}(u)\ldots\ge J_{n-1}^{(n-1)}(u)$ are distributed like the decreasing order statistics of $n-1$ independent and identically distributed  random variables $(J_i(u))_{1\le i\le n-1}$, each having the distribution in \eqref{Jdis0}.

 {\rm (ii)}\ 
For $n$, $r$, $x_k$, $z$ as specified, 
\bea\label{top-}
&&
\lim_{t\dto 0}
\PP\bigg(\frac{\Delta_t^{(r+k)} }{\Delta_t^{(r+n)}}>x_k,\ 0\le k\le n-1\bigg|
\Delta_t^{(r+n)}=\pibarinv(z/t)\bigg)\cr
&&=\
\PP\big(J_{n-1}^{(k)}\big({\rm B}^{1/\alpha}_{r,n}\big)>x_k, 1\le k\le n-1,\, {\rm B}^{1/\alpha}_{r,n}\le x_0^{-1}\big),
\eea
where the 
$J_i(u)$ are as in \eqref{13b} and $B_{r,n}$ is a Beta$(r,n)$  random variable  independent of  
$(J_i(u))_{1\le i\le n-1}$.

 {\rm (b)}\
When $\alpha=0$, each ratio $\Delta^{(r+k)}_t/\Delta_t^{(r+n)}\topr \infty$ as $t\dto 0$,
for $1\le k\le n-1$.
When $\alpha=\infty$, each ratio  $\Delta^{(r+k)}_t/\Delta_t^{(r+n)}\topr 1$ as $t\dto 0$,
for $1\le k\le n-1$.
\end{proposition}

\medskip\noindent{\bf Proof of Proposition \ref{cor1}:}\ 
Equations \eqref{13b} and \eqref{top-} are implicitly proved in the proof of Theorem \ref{preT_2}.
Part (b) follows from similar arguments as in Part (ii) of Proposition \ref{cor2a}. \halmos

\begin{remark}\label{rem6} {\rm
(i)\ In Part  (a)(i) of Proposition \ref{cor1} the $x_0$ variable is superfluous, but it is relevant in Part (a)(ii).

(ii)\ If we make the convention that $B_{0,n}\equiv 0$ a.s., put $u=0$ in \eqref{Jdis0}, and identify $(J_i(0))$ with  a sequence $(L_i)$ of
independent and identically distributed random variables each having the  distribution  
defined in \eqref{Ldis},
 we get the case $r=0$ of \eqref{top-}; namely, 
for  $x_k\ge 1$, $0\le k\le n-1$, $n=2,3,\ldots$, and $z>0$,
\be\label{top+}
\PP\left(\frac{\Delta_t^{(k)} }{\Delta_t^{(n)}}>x_k, 1\le k\le n-1\bigg|
\Delta_t^{(n)}=\pibarinv(z/t)\right)
\to
\PP\big(L_{n-1}^{(k)}>x_k, 1\le k\le n-1\big),
\ee
as $t\dto 0$, 
where 
$L_{n-1}^{(1)}\ge L_{n-1}^{(2)}\ldots\ge L_{n-1}^{(n-1)}$ are  the decreasing order statistics of $(L_i)_{1\le i\le n-1}$.
\sid{Equation} \eqref{top+} can of course be proved directly. 

(iii)\ 
The case $r\in\N$, $n=1$, in Part(a) (i) of Proposition \ref{cor1},
is covered by setting $n=r+1$, and $x_1=\cdots=x_{r-1}=1$ when $r>1$, in  \eqref{top+}, to get
\be\label{n=1}
\lim_{t\dto 0}
\PP\left(\frac{\Delta_t^{(r)} }{\Delta_t^{(r+1)}}>x_r\bigg|
\Delta_t^{(r+1)}=\pibarinv(z/t)\right)=\PP\big(L_r^{(r)}>x_r\big)=
x_r^{-r\alpha},
\ee
for $x_r\ge 1$ and  $z>0$.
Here $L_r^{(r)}\eqdr \min_{1\le i\le r}L_i$, where  $(L_i)_{1\le i\le r}$ are i.i.d.  random variables, each having the  distribution in \eqref{Ldis}. 
Note that $L_r^{(r)}\eqdr {\rm B}_{r,1}^{-1/\alpha}$.

(iv)\ Convergence of the conditional distributions in \eqref{top-}, \eqref{top+}, and \eqref{n=1}, together with
\be\label{nc}
\lim_{t\dto 0} \PP\big(
t\pibar(\Delta^{(r+j)}_t)\le z_{r+j},\, 0\le j\le n\big)
=\PP\left(\Gamma_{r+j}\le z_{r+j},\, 0\le j\le n\right),
\ee
 for $0\le z_{r}\le \cdots \le z_{r+n}$,
 implies convergence of the corresponding joint, and hence marginal, distributions.
Since the right-hand sides of  \eqref{top-}, \eqref{top+}, and \eqref{n=1} do not depend on $z$, independence obtains in the corresponding limiting joint distributions.
To verify \eqref{nc}, observe
 we have for $j = 0, 1,\ldots, n$, 
\begin{align*}
t\pibar(\Delta_t^{(r+j)}) \eqdr t\pibar(\pibarinv(\Gamma_{r+j}/t))
\sim \frac{\pibar\big(\pibarinv\big(\frac{\Gamma_{r+j}}{t} \big)\big)}{\pibar\big(\pibarinv\big(\frac{1}{t}\big)\big)}
\sim \frac{\Gamma_{r+j}/t}{1/t}  = \Gamma_{r+j},
\end{align*}
where  the convergence is almost sure
as $t\dto 0$.
}\end{remark}

\begin{proposition}[Ratios smaller than  1]\label{preT2} 
Suppose $\pibar(\cdot) \in RV_0(-\alpha)$ with $0<\alpha<\infty$ and $r,n \in\N$. 

{\rm (i)}\ For each $z>0$ and $w\in(0,1)$ 
\bea\label{13c}
&&
\lim_{t\dto 0}
\PP\Big( W_{r,n}(t)\le w \, \big| \, \Delta_t^{(r+n)}= \pibarinv(z/t) \Big)
= \PP \big(K_{r+n-1}^{(n)}\le w \big),
\eea
where
$K_{r+n-1}^{(n)}$ is the $n$th largest 
of $r+n-1$ independent and identically distributed  random variables $(K_i)_{1\le i\le r+n-1}$, with distribution 
$ \PP(K_1 \le w) = w^{\alpha}$,  $w \in (0,1)$.

{\rm (ii)}\  For each $z>0$ and $w\in(0,1)$ 
\be\label{cdu}
\lim_{t\dto 0}
\PP\left(\Delta_t^{(r)}\ge \pibarinv(z/t) \, \big| \,W_{r,n}(t)=w\right)
= \PP \left(\Gamma_{r+n} \le w^{-\alpha}z\right).\ee

\end{proposition}

\begin{remark}\label{Brn_dis} {\rm 
(i)\ 
Taking expectations in \eqref{13c} gives, as $t\downarrow 0$,
\ben  
W_{r,n}(t)=
\frac{ \Delta_t^{(r+n)}}{\Delta_t^{(r)}}
\todr 
 K_{r+n-1}^{(n)}.
\een
The connection with  \eqref{Ldis} is that $L \eqd 1/K$. 
From \eqref{top+} with $n$ replaced by $r+n$, 
 we thus have the various alternatives 
\ben 
W_{r,n}(t)\todr 
W_{r,n} \eqd B_{r,n}^{1/\alpha} \eqd K_{r+n-1}^{(n)} \eqd 1/L_{r+n-1}^{(r)}, \ {\rm as}\ t\dto 0,
\een
where $L_{r+n-1}^{(r)}$ is the $r$th largest 
of  i.i.d. rvs
$(L_i)_{1\le i \le r+n-1}$.

}
\end{remark}

\medskip\noindent {\bf Proof of  Proposition \ref{preT2}:}\
(i)\ 
 The probability on the LHS of \eqref{13c} equals
\ben
\PP\big(\Delta_t^{(r)}\ge w^{-1}\pibarinv(z/t)\big|\Delta_t^{(r+n)}= \pibarinv(z/t)\big).
\een
Conditional on $\{\Delta_t^{(r+n)}= \pibarinv(z/t)\}$, the ordered
$\Delta_t^{(1)}\ge \cdots\ge \Delta_t^{(r+n-1)}$
have the distribution of the decreasing order statistics 
$K_{r+n-1}^{(1)}(t,z)\ge\cdots\ge K_{r+n-1}^{(r+n-1)}(t,z)$  of $r+n-1$ independent and identically distributed random variables $(K_i(t,z))_{i=1,2,\ldots}$, each having the distribution 
\ben 
\PP(K_1(t,z)\in\rmd x)= 
\frac{\Pi(\rmd x)
{\bf 1}\{x \ge \pibarinv(z/t)\}}
{\pibar(\pibarinv(z/t)-)},\ x>0.
\een
From  \eqref{1b} it follows that, as $t\dto 0$, for each $w\in(0,1)$ and $z>0$,
\bean 
\PP\left(K_1(t,z) > w^{-1}\pibarinv(z/t)\right)
&=&
\frac{t\pibar(w^{-1}\pibarinv(z/t)-)
{\bf 1}\{w^{-1}>1\}}
{t\pibar(\pibarinv(z/t)-)}  \cr
&&\cr
&\to &
\frac{w^{\alpha}z{\bf 1}_{\{w<1\}}}{z}\cr
&=&
w^{\alpha} {\bf 1}_{\{w<1\}}=\PP(K_1\le w),
\eean
where $K_1$ is a  random variable such that $K_1^\alpha$ is $U[0,1]$.
Thus as $t\dto 0$, for $0<w<1$,
\bean 
&&\PP\Big(\Delta_t^{(r)}\ge w^{-1}\pibarinv(z/t)\, \Big| \,\Delta_t^{(r+n)}= \pibarinv(z/t)  \Big)\cr
&&\cr
&&= 
\PP\left(\, {\rm at\ least}\ r\ {\rm of}\ K_1(t,z), \ldots, K_{r+n-1}(t,z)\ \text{ exceed or equal}\ w^{-1}\pibarinv(z/t) \, \right)\cr
&&\cr
&&=
\sum_{k=r}^{r+n-1} {r+n-1\choose k}
\left(\frac{\pibar(w^{-1}\pibarinv(z/t)-)}  {\pibar(\pibarinv(z/t)-)} \right)^k
\left(1-\frac{\pibar(w^{-1}\pibarinv(z/t)-)}{\pibar(\pibarinv(z/t)-)} \right)^{r+n-1-k}\cr
&&\cr
&&\to
\sum_{k=r}^{r+n-1} {r+n-1\choose k}
\left(w^\alpha\right)^k
\left(1-w^\alpha \right)^{r+n-1-k}\cr
&&\cr
&&= 
\PP\left({\rm at\ least}\ r\ {\rm of}\ K_1,\ldots, K_{r+n-1}\ {\rm are\ smaller\ than}\ w\right),
\eean
and this is the RHS of \eqref{13c}.  

(ii)\ 
Using \eqref{ch} suggests writing  the LHS of \eqref{cdu} as
\begin{align*}
\PP\Big( t\pibar\big(\pibarinv(\Gamma_r/t)\big) \le z
\, \big|\,  t\pibar\big(\pibarinv(\Gamma_r + \wt \Gamma_n)/t\big) 
= t\pibar\big(w\pibarinv(\Gamma_r/t)\big)\Big),
\end{align*}
where $\wt\Gamma_n$ is a Gamma$(n,1)$ rv independent of $\Gamma_r$, 
a Gamma$(r,1)$ rv.
By \eqref{1b} it is plausible  that this tends to 
$ \PP\big(\Gamma_r\le z \, \big| \,  \Gamma_r+\wt\Gamma_n= w^{-\alpha}\Gamma_r\big)$ as $t\dto 0$. To prove it, 
write, to be brief, $A_t:= t\pibar(\pibarinv(\Gamma_r/t))$, $W_t:=W_{r,n}(t)$ and $W:=W_{r,n}$, 
 and let $C$ be any Borel subset of $\R^+$. Then, for each $t>0$,  by dominated convergence,
\bean
&&
\int_{w\in C} \PP(A_t\le z|W_t=w)\PP(W\in \rmd w)\cr
&&\cr
&&=\int_{w\in C}\lim_{\veps\dto 0}\left(
\frac{\PP(A_t\le z, w-\veps<W_t<w+\veps)}{\PP(w-\veps<W_t<w+\veps)}\right) \PP(W\in \rmd w)\cr
&&\cr
&&=\lim_{\veps\dto 0}\int_{w\in C}\left(
\frac{\PP(A_t\le z, w-\veps<W_t<w)}{\PP(w-\veps<W_t<w)}\right) \PP(W\in \rmd w).
\eean
Since $(A_t,W_t)$ converges in distribution to $ (\Gamma_r,W)$, as $t\dto 0$, and the limit distribution is continuous, the convergence is uniform. So, given $\delta>0$, there is a $t_0=t_0(\veps, \delta)>0$ such that the last expression is, for $0<t\le t_0$, no greater than 
\bean
&&\limsup_{\veps\dto 0}\int_{w\in C}\left(
\frac{\PP(\Gamma_r\le z, w-\veps<W<w+\veps)}{\PP(w-\veps<W<w+\veps)}+\delta\right) \PP(W\in \rmd w)\cr
&&\cr
&&\le \int_{w\in C}
\PP(\Gamma_r\le z| W=w) \PP(W\in \rmd w)+\delta.
\eean
In a similar way we find a lower bound for the liminf with $-\delta$, hence, for $0<t\le t_0$,
\ben
\left|\int_{w\in C}\big(
\PP(A_t\le z|W_t=w)-
\PP(\Gamma_r\le z| W=w)
\big)\PP(W\in \rmd w)\right|\le \delta.
\een
Now, choosing Borel sets $C^+$ and $C^-$ on which the integrand in the last integral is positive or negative, we see that
\be\label{d2}
\int_{w>0}\left|
\PP(A_t\le z|W_t=w)-
\PP(\Gamma_r\le z| W=w)\right|
\PP(W\in \rmd w)\le \delta
\ee
for $0<t\le t_0$. 
Take any sequence $t_k\dto 0$.
Let $k\to\infty$ and  use Fatou's lemma, then let $\delta\dto 0$, to deduce from 
\eqref{d2} that 
\ben
\liminf_{k\to\infty}\PP(A_{t_k}\le z|W_{t_k}=w)=
\PP(\Gamma_r\le z| W=w)
\een
(a.e. $w$ with respect to the distribution of $W$, thus, Lebesque a.e.).
Taking a further subsequence if necessary, we can replace ``liminf" by ``lim" here.
Then, since the limit holds for arbitrary $t_k$, we conclude that
\ben 
\lim_{t\dto 0}\PP(A_t\le z|W_t=w)=
\PP(\Gamma_r\le z| W=w), \ {\rm a.e.}\ (w).
\een

We can evaluate the probability on the RHS here
using $W\eqdr \{(\Gamma_r+\wt\Gamma_n)/\Gamma_r\}^{-1/\alpha}$ and 
 \ben
  \PP\big(\Gamma_r\le z,\,  \Gamma_r+\wt\Gamma_n\ge w^{-\alpha}\Gamma_r\big)=
 \int_{0\le v\le z} \int_{x \ge w^{-\alpha }v-v}
 \left( \frac{ e^{-x}x^{n-1}}  { \Gamma(n)}\right) \rmd x 
 \left(\frac{e^{-v}v^{r-1}}{\Gamma(r)}\right)\rmd v.
\een
Differentiate with respect to $w$  and divide by the density of $W_{r,n}$  to get the required limiting conditional distribution as
\ben
\frac{1}{f_{W_{r,n}}(w)}
\int_{v=0}^z
 \left(  \frac{e^{-(w^{-\alpha}-1)v}(w^{-\alpha }-1)^{n-1}\alpha w^{-\alpha-1} v^n}
 {\Gamma(n)}\right)
 \left(\frac{ e^{-v}v^{r-1}}{\Gamma(r)}\right)\rmd v.
 \een
 Substituting for $f_{W_{r,n}}(w)$ from \eqref{prodBeta2}, we can calculate 
the last expression as 
\[
\int_0^z w^{-\alpha(r+n)} \frac{e^{-w^{-\alpha}v} v^{r+n-1}}{\Gamma(r+n)} \rmd v = \int_0^{w^{-\alpha}z} e^{-v} \frac{v^{r+n-1}}{\Gamma(r+n)} \rmd v,
\]
which is the RHS of  \eqref{cdu}. 
\halmos

\begin{proposition}[Ratios greater than  1]\label{sum_l}
Let \sid{$\{\Gamma_j^{-1/\alpha}, j\geq 1\}$}
be the ordered points of a Poisson  point process with intensity measure $\Lambda(\rmd x) = \alpha x^{-\alpha-1}\rmd x {\bf 1}_{x > 0}$, $\alpha>0$.
 Let $ r, n \in\N$, $0<u<1$, $\lambda>0$. 
 Then we have the conditional Laplace transform 
\be\label{iii}
\EE\left(\exp\Biggl(-\lambda  \sum_{i=1}^{n-1} \sid{
    \Bigl(\frac{\Gamma_{r+i}}{ \Gamma_{r+n}}\Bigr)^{-1/\alpha}} \right)
\, \bigg| \,\sid{ \Bigl(\frac{\Gamma_{r+n}}{\Gamma_r }  \Bigr)^{-1/\alpha}}
= u\Biggr) =(\Phi(\lambda,u))^{n-1}, 
\ee
where 
\[
\Phi(\lambda,u)= \frac{\int_{1}^{1/u}    e^{-\lambda x} \Lambda(\rmd x)}{1- u^{\alpha}}.
\] 

When $r=0$, the sum of ratios of jumps greater than 1 has representation 
\be\label{D2}
\sid{\sum_{i=1}^{n-1} \Bigl(\frac{\Gamma_{i}}{ \Gamma_{n}}\Bigr)^{-1/\alpha} }\eqdr \sum_{i=1}^{n-1} L_i,
\ee
where  
the $L_i$ are i.i.d random variables each with the same distribution
as 
\sid{$\Gamma_{1}/\Gamma_{2}$;} namely, $\PP(L_1\in\rmd x)= \Lambda(\rmd x) {\bf 1}_{\{x>1\}}$. 
\end{proposition}

\begin{remark}{\rm
The representation \eqref{D2} of the sum of ratios of the ordered jumps of a stable subordinator as a random walk in $n$ provides the impetus for further work in large trimming results in the spirit of the investigations in \cite*{BMR2016}. 
}\end{remark}

\medskip\noindent {\bf Proof of  Proposition \ref{sum_l}:}\
Equality \eqref{iii} can be read from \res{the order statistics
  property of the homogeneous Poisson process} when $r > 0$ and from \eqref{Mt1} when $r = 0$.

\halmos


\section{Converse Results}\label{s5}
 Theorem \ref{conv} gives converses to the previous results.

\begin{theorem}[Converse Results: Ratios Bigger than 1.]\label{conv}
\noindent 
Suppose, for some $r\in\N$, $n\in\N$,   $\Delta^{(r)}_t/\Delta_t^{(r+n)}\todr Y$, as
$ t\dto 0$, for an extended value random variable\footnote{\
A random variable that may take the value $+\infty$ with positive probability.}
 $Y\ge 1$. Then one of the following holds:

\ {\rm (i)}\ $\PP(1<Y<\infty)>0$, in which case 
$\pibar(\cdot) \in RV_0(-\alpha)$ with $0<\alpha<\infty$; 

\ {\rm (ii)}\   $\PP(Y=1)=1$, in which case $\pibar$ is rapidly varying at 0;

\ {\rm (iii)}\ $Y=\infty$ a.s., in which case $\pibar$ is slowly varying at 0.

\begin{remark}{\rm [Converse Results: Ratios Smaller than 1.]
Analogous results to  Theorem \ref{conv} for ratios smaller than 1
follow by taking reciprocals.
Write $\Delta^{(r+n)}_t/\Delta^{(r)}_t=
(\Delta^{(r)}_t/\Delta^{(r+n)}_t)^{-1}$
and apply the theorem, replacing $Y$ by $Y^{-1}$, and making the obvious interpretations in Parts (i), (ii) and (iii) of the theorem.
}
\end{remark}
%
%
%
\end{theorem}

\medskip\noindent {\bf Proof of  Theorem \ref{conv}:}\
Assume for some $r\in\N$,  $n\in\N$, 
\be\label{C1a}
\frac{\Delta^{(r)}_t}{\Delta^{(r+n)}_t}
\todr Y, \ {\rm as}\ t\dto 0,
\ee
where $Y$ is an extended  random variable with distribution $G$, say, on $[1,\infty]$. 
The proof that follows is similar in style to that of \cite{KeveiMason2014} who considered ratios of
successive  jumps, that is, the case $n=1$. When $n >1$, some rather different arguments are needed at some places.

Keep $u$ fixed in $(0,1)$ throughout the remainder of the proof and 
use \eqref{ch} to write  
\begin{align*}  
\PP\Big(\Delta^{(r+n)}_t <  u\Delta_t^{(r)}\Big)&=
\PP\Big(\pibarinv((\Gamma_r+ \wt\Gamma_{n} )/t)< u\pibarinv(\Gamma_r/t) \Big) \cr
&=
\int_{y\ge 0}
\PP\Big(\wt \Gamma_{n} >  t\pibar\left(u\pibarinv(y/t)\right)-y\Big)\PP \big(\Gamma_r\in\rmd y\big),
\end{align*}
where $\Gamma_{r}$ and $\wt \Gamma_n$ are independent Gamma  random variables.
Substituting for their densities gives
\bea\label{C3}
&&
\PP\Big(\Delta^{(r+n)}_t < u\Delta_t^{(r)}\Big)=
\int_{y >  0}\int_{z >  t\pibar\left(u\pibarinv(y/t)\right)-y}
\left(\frac{e^{-z}z^{n-1}}{\Gamma(n)}\right)\rmd z 
\left(\frac{e^{-y}y^{r-1}}{\Gamma(r)}\right)\rmd y
\cr
&&\cr
&&\cr
&&=
\frac{t^{r+n}}{\Gamma(r+n)}
\int_{y> 0}\int_{z > \pibar\left(u\pibarinv(y)\right)}
\left(\frac{e^{-tz}(z-y)^{n-1}}{B(r,n)}\right)\rmd z \,  y^{r-1}\rmd y  \cr
&&\cr
&&\cr
&&=
\frac{t^{r+n}}{\Gamma(r+n)}
\int_{z> 0}\int_{y < \pibar\left(\pibarinv(z)/u\right)/z} 
\left(\frac{(1-y)^{n-1} y^{r-1}}{B(r,n)}\right)\rmd y\, 
e^{-tz}z^{r+n-1}\rmd z.
\eea
Here note that, since $u<1$, we have
$t\pibar\left(u\pibarinv(y/t)\right)\ge 
t\pibar\left(\pibarinv(y/t)-\right)\ge y$, 
and
$\pibar\left(\pibarinv(z)/u\right)/z\le 
\pibar\left(\pibarinv(z)\right)/z\le 1$.
We recognise the inner integral in \eqref{C3} as the incomplete Beta function
$B\left(r,n; \pibar\left(\pibarinv(z)/u\right)/z\right)$ (see \eqref{incombet}).
By assumption \eqref{C1a}, the expression in \eqref{C3} tends to $\Gbar(1/u):=P(Y > 1/u)$
as $t\dto 0$,  at continuity points of $\Gbar$. To simplify the notation, from this point on let $x=1/u>1$.
Let 
\ben
U_x(z):= \int_0^z B\left(r,n; \pibar\left(\pibarinv(v)x\right)/v\right)v^{r+n-1}\rmd v,\ z>0.
\een
Then from \eqref{C3}, 
\ben
t^{r+n}\int_{z>0}e^{-tz}U_x(\rmd z)\to
\Gamma(r+n)\Gbar(x),\ {\rm as}\ t\dto 0,
\een
 at continuity points of $\Gbar$,
which by Thm. 1.7.1 p.37 of \cite{BGT87} implies
\ben
z^{-r-n}U_x(z) \to \Gamma(r+n)\Gbar(x)/\Gamma(r+n+1)
=\Gbar(x)/(r+n),\ {\rm as}\ z\to\infty.
\een
Write this as
\be\label{C4}
\frac{1}{z^{r+n}} \int_0^z b_x(v)v^{r+n-1}\rmd v
\to \frac{\Gbar(x)}{r+n},\ {\rm as}\ z\to\infty,
\ee
where 
\ben
b_x(v)=B(r,n;f_x(v))=
\frac{1}{B(r,n)}\int_0^{f_x(v)}y^{r-1}(1-y)^{n-1}\rmd y
=: \int_0^{f_x(v)}p(y)\rmd y,
\een
with 
$f_x(v):= \pibar\left(\pibarinv(v) x\right)/v$, $v>0$ 
and $p(y):= y^{r-1}(1-y)^{n-1}/B(r,n)$, $0\le y\le 1$.
Note that $x$ is kept fixed in $b_x(v)$ and $f_x(v)$.
We have $0\le f_x(v)\le 1$, so $0\le b_x(v)\le 1$, for all $v>0$.
\eqref{C4} implies
\be\label{C5}
\frac{1}{z^{r+n}} \int_z^{\lambda z} b_x(v)v^{r+n-1}\rmd v
=\int_1^\lambda b_x(vz)v^{r+n-1}\rmd v
\to \frac{(\lambda^{r+n}-1)\Gbar(x)}{r+n},\ {\rm as}\ z\to\infty,
\ee
for any $\lambda>1$ and each fixed $x>0$. 

Functions $b_x(v)$, $f_x(v)$, are not necessarily monotone but are of bounded variation (BV) on finite intervals bounded away from 0.
To see this, observe that the function
$m_x(v):= vf_x(v)=\pibar\left(\pibarinv(v)x\right)$ is nondecreasing in $v$ and
\ben
|\rmd f_x(v)|=\bigg|\frac{\rmd m_x(v)}{v}-\frac{m_x(v)\rmd v}{v^2}\bigg|
\le \frac{\rmd m_x(v)}{v}+\frac{\rmd v}{v};
\een
thus, with $p_0:=\sup_{0\le y\le 1}p(y)$, 
\ben
|\rmd b_x(v)|=|p(f_x(v))\rmd f_x(v)|\le
p_0\left( \frac{\rmd m_x(v)}{v}+\frac{\rmd v}{v}\right),
\een
and the RHS is integrable over $v\in[\delta,z]$, for any $0<\delta<z$. So $f_x$ and $b_x$ are of bounded variation on $[\delta,z]$ for any $0<\delta<z$.
Take any sequence $z_k\to\infty$. By Helly's theorem for finite measures 
we can find a subsequence, also denoted $z_k$, possibly depending on $x$, such that
\ben
b_x(vz_k)\to g_x(v),\  v>0, \ {\rm as}\ k\to\infty,
\een
at continuity points of $g$,
for a function $g_x(v)\in[0,1]$.
Using dominated convergence in \eqref{C5} we get
\ben 
\int_1^\lambda g_x(v)v^{r+n-1}\rmd v
= \frac{(\lambda^{r+n}-1)\Gbar(x)}{r+n}=\Gbar(x)
\int_1^\lambda v^{r+n-1}\rmd v.
\een
This holds for all $\lambda>1$ and so implies $g_x(v)=\Gbar(x)$, for all $v>1$, $x>0$, not depending on the choice of subsequence. Thus we deduce that
\ben 
b_x(vz)=\int_0^{f_x(vz)}p(y)\rmd y\to \Gbar(x),
\een
as $z\to\infty$, at continuity points of $\Gbar$, for all $v>1$.
Take $v=2$. 
Now $f_x(2z)$ is monotone in $x$ for each $z$, 
so by Helly's theorem again each sequence $z_k\to\infty$ contains a further subsequence,  also denoted $z_k$,
such that $f_x(2z_k)\to h(x)\in[0,1]$, as $k\to\infty$,  at continuity points of $h(x)$.
Thus we obtain
\be\label{hdef}
\int_0^{h(x)} p(y)\rmd y=\Gbar(x)
\ee
at continuity points of $h$. Again the limit does not depend on the choice of subsequence.
This identifies $h(x)$ as $\Iinv(\Gbar(x))$, where $\Iinv(\cdot)$ is the unique inverse function to the continuous strictly increasing function
$I(\cdot)=\int_0^\cdot p(y)\rmd y$.
Thus, continuity points of $h$ are points of increase of $G$.
Define 
\[
\AAA := \{ x \ge 1: x \text{ is a continuity point and a point of increase of } G \}.
\]
We conclude that
\be\label{C8}
 \lim_{z\to\infty}\frac{\pibar\left(\pibarinv(z) x\right)}{z}
 =\lim_{z\to\infty}f_x(z)= \lim_{z\to\infty}f_x(2z)=h(x),
 \ {\rm for\ all}\ x \in \AAA,
\ee
where $h$ satisfies \eqref{hdef}.
\eqref{C8} is exactly analogous to Eq.(2.10) of \cite{KeveiMason2014} and we follow their arguments henceforth to finish the converse part of the proof. There are three alternatives.

(i)\ $\PP(1<Y<\infty)>0$.
In this  case $\Gbar$ has at least one point of decrease in $(1,\infty)$, say $x$, and a neighbourhood $(x-\veps, x+ \veps)$ for some $\veps> 0$, such that $\Gbar(y) > 0 $ for all $y$ in the neighbourhood. \cite{KeveiMason2014} gave a careful analysis of this situation, showing that it leads to $\pibar(\cdot) \in RV_0(-\alpha)$ with $0<\alpha<\infty$.

\medskip
(ii)\   $\PP(Y=1)=1$. This means that $\PP(Y>x)=\Gbar(x)=0$ for all $x>1$, so $\int_0^{h(x)} p(y)\rmd y=0$ and 
\ben
 \lim_{z\to\infty}\frac{\pibar\left(\pibarinv(z)x\right)}{z}\le 
\lim_{z\to\infty}\frac{\pibar\left(\pibarinv(z)x\right)}{\pibar\left(\pibarinv(z)\right)}
 =0
\een
 for all $x>1$. Thus $\pibar$ is rapidly varying at 0.

\medskip
(iii)\  $Y=\infty$ a.s.
This means that $\PP(Y>x)=\Gbar(x)=1$ for all $x>1$, so $\int_0^{h(x)} p(y)\rmd y=1$ and 
\ben
 \lim_{z\to\infty}\frac{\pibar\left(\pibarinv(z)x\right)}{z}=1
\een
for all $x>1$. This leads to $\pibar$ slowly varying at 0 as shown in \cite{KeveiMason2014}, and completes the proof.
 \halmos

\section{Related Results: Order Statistics of i.i.d. rvs}\label{s6}
  We conclude with some history relating how these kinds of results have antecedents in the  literature of order statistics 
 of i.i.d. real-valued random variables. 
 The general scenario there is of 
 the order statistics $X_n^{(n)} \le \cdots \le  X_n^{(1)}$ of i.i.d.  rvs $(X_i)_{i\in\N}$ in $\R$ with distribution $F$ such that $F(x)<1$ for all $x$. The asymptotic is then as $n\to \infty$ (``large time").
(In most of the results quoted below  the distribution $F$ is also assumed  continuous, so ties among order statistics have probability 0. We  avoided such an assumption on $\Pi$ in our results.)
 
 An early and well-cited venture in this area 
 was by 
 \cite{arov1960}.   
 They considered not only the order statistics but also their sum, i.e., the random walk $S_n$ whose step sizes are the $X_i$, obtaining among other things results for convergence of joint distributions of deterministically normed order statistics, and as a corollary limiting distributions for ratios of (not necessarily successive) order statistics.  
 This was extended to ratios of the sum after removal of a fixed number of extreme terms (i.e., the trimmed sum) to large order statistics. 
 (The distribution $F$ was assumed to have a density.)

 \cite{smidstam1975}
considered the $(X_n^{(i)})_{1\le i\le n}$ as above, 
and, in what amounts to a generalisation of and converse to one of the  \cite{arov1960} results, showed that 
$\lim_{n\to\infty}
\PP(X_n^{(j+1)}/X_n^{(j)}\le x,\,
1\le j\le k) =\prod_{j=1}^k x^{j\alpha}$
for all $x\in (0,1)$ and $k\in\N$
iff
the distribution tail 
$\Fbar(x)\in RV_\infty(-\alpha)$, $\alpha\ge 0$. 
They include the $\alpha=0$ case ($\Fbar$ slowly varying at $\infty$). 
Their proof used Scheff\'e's lemma 
and applications of the Wiener-Tauberian theory. 
     A converse to another of the  \cite{arov1960}   results is in \cite{malres1984}.
An earlier result along the lines of \cite{smidstam1975}
is in \cite{shorrock1972}.  

\cite{teugels1981} considered order statistics of i.i.d. rvs in the domain of attraction of a stable law, and gave results extending some of the Arov-Bobrov limit laws concerning ratios of sums of order statistics to their (trimmed) sums. 
For an application of these kinds of ideas in
 reinsurance, see  \cite{ladteu2006}.  

\cite{lanstad2002}  
gave a simplified version of the 
 \cite{smidstam1975}  
 result (for the $k=1$ case)
and  extended this for when $F$ is in the domain of attraction of an extreme value distribution. 

There is of course in addition a very large literature analysing various functions of order statistics  of i.i.d. real-valued rvs which we do not attempt to summarise here.

We remark finally that while there are obvious correspondences between the (large-time) 
i.i.d. case and the (small time) point process case, there are significant differences too.
One aspect is that, in view of our assumption  $\pibar(0+)=\infty$, there are always infinitely many points of the process  in any right neighbourhood of $0$, hence,  infinitely many ordered  points; whereas, in the i.i.d. case, there are of course at most $n$ order statistics in a sample of size $n$.
Thus there is no immediate counterpart of results like \eqref{Mt0} or \eqref{asymIndep}.
This feature actually simplifies some of the point process proofs, for example that  of  Theorem \ref{preT_2}, although the formulation is more complex. 

\medskip\noindent {\bf Acknowledgements.}\
We are grateful for helpful feedback from 
P\'eter Kevei and David Mason.

\renewcommand{\bibfont}{\small}
\bibliography{Library_Levy_July2017}

\begin{thebibliography}{}

\bibitem[\protect\citeauthoryear{Arov \& Bobrov}{Arov \&
  Bobrov}{1960}]{arov1960}
Arov, D. \& Bobrov, A. (1960).
\newblock The extreme terms of a sample and their role in the sum of
  independent variables.
\newblock {\em Theory Probab. Appl.}, {\em 5}, 377--396.

\bibitem[\protect\citeauthoryear{Bertoin}{Bertoin}{2006}]{Bertoin2006}
Bertoin, J. (2006).
\newblock {\em Random Fragmentation and Coagulation Processes}.
\newblock Cambridge studeis in advanced mathematics 102, Cambridge University
  Press, Cambridge.

\bibitem[\protect\citeauthoryear{Bingham, Goldie \& Teugels}{Bingham
  et~al.}{1987}]{BGT87}
Bingham, N.~H., Goldie, C.~M., \& Teugels, J.~L. (1987).
\newblock {\em Regular Variation}.
\newblock Cambridge University Press.

\bibitem[\protect\citeauthoryear{Buchmann, Fan \& Maller}{Buchmann
  et~al.}{2016}]{bfm2016}
Buchmann, B., Fan, Y., \& Maller, R.~A. (2016).
\newblock Distributional representations and dominance of a {L}\'evy process
  over its maximal jump processes.
\newblock {\em Bernoulli}, {\em 22\/}(4), 2325--2371.

\bibitem[\protect\citeauthoryear{Buchmann, Maller \& Resnick}{Buchmann
  et~al.}{2016}]{BMR2016}
Buchmann, B., Maller, R.~A., \& Resnick, S.~I. (2016).
\newblock Processes of rth largest.
\newblock {\em arXiv:1607.08674}.

\bibitem[\protect\citeauthoryear{Ferguson \& Klass}{Ferguson \&
  Klass}{1972}]{ferguson:klass:1972}
Ferguson, T. \& Klass, M. (1972).
\newblock A representation of independent increment processes without gaussian
  components.
\newblock {\em Ann. Math. Statist..}, {\em 43\/}(5), 1634--1643.

\bibitem[\protect\citeauthoryear{Gregoire}{Gregoire}{1984}]{Gregoire1984}
Gregoire, G. (1984).
\newblock Negative binomial distributions for point processes.
\newblock {\em Stochastic Process. Appl.}, {\em 16\/}(2), 179--188.

\bibitem[\protect\citeauthoryear{Ipsen \& Maller}{Ipsen \&
  Maller}{2017a}]{stableJump1}
Ipsen, Y.~F. \& Maller, R.~A. (2017a).
\newblock Convergence to stable limits for ratios of trimmed {L}\'evy processes
  and their jumps.
\newblock {\em unpublished manuscript}.

\bibitem[\protect\citeauthoryear{Ipsen \& Maller}{Ipsen \&
  Maller}{2017b}]{trimfrenzy}
Ipsen, Y.~F. \& Maller, R.~A. (2017b).
\newblock Generalised {P}oisson-{D}irichlet distributions and the negative
  binomial point process.
\newblock {\em arXiv:1611.09980}.

\bibitem[\protect\citeauthoryear{Kevei \& Mason}{Kevei \&
  Mason}{2014}]{KeveiMason2014}
Kevei, P. \& Mason, D.~M. (2014).
\newblock The limit distribution of ratios of jumps and sums of jumps of
  subordinators.
\newblock {\em Lat. Am. J. Probab. Math. Stat.}, {\em 11\/}(2), 631--642.

\bibitem[\protect\citeauthoryear{Kingman}{Kingman}{1975}]{Kingman1975}
Kingman, J. F.~C. (1975).
\newblock Random discrete distributions.
\newblock {\em J. R. Stat. Soc. Series B Stat. Methodol.}, {\em 37\/}(1),
  1--22.

\bibitem[\protect\citeauthoryear{Ladoucette \& Teugels}{Ladoucette \&
  Teugels}{2006}]{ladteu2006}
Ladoucette, S.~A. \& Teugels, J.~L. (2006).
\newblock Reinsurance of large claims.
\newblock {\em J. Comput. Appl. Math.}, {\em 186\/}(1), 163--190.

\bibitem[\protect\citeauthoryear{Lanzinger \& Stadtm\"{u}ller}{Lanzinger \&
  Stadtm\"{u}ller}{2002}]{lanstad2002}
Lanzinger, H. \& Stadtm\"{u}ller, U. (2002).
\newblock Tauberian theorems and limit distributions for upper order
  statistics.
\newblock {\em Publ. Inst. Math. Nouvelle Ser.}, {\em 71}, 41--53.

\bibitem[\protect\citeauthoryear{LePage}{LePage}{1980}]{LeP1980a}
LePage, R. (1980).
\newblock Multidimensional infinitely divisible variables and processes {P}art
  {I}.
\newblock {\em Technical Rept. 292 Dept. Statistics, Stanford University}.

\bibitem[\protect\citeauthoryear{LePage}{LePage}{1981}]{LeP1981}
LePage, R. (1981).
\newblock Multidimensional infinitely divisible variables and processes {P}art
  {II}.
\newblock In {\em Probability in Banach Spaces {III}}  (pp.\ 279--284).
  Springer.

\bibitem[\protect\citeauthoryear{LePage, Woodroofe \& Zinn}{LePage
  et~al.}{1981}]{LePZinn1981}
LePage, R., Woodroofe, M., \& Zinn, J. (1981).
\newblock Convergence to a stable distribution via order statistics.
\newblock {\em Ann. Probab.}, {\em 9}, 624--632.

\bibitem[\protect\citeauthoryear{Maller \& Resnick}{Maller \&
  Resnick}{1984}]{malres1984}
Maller, R. \& Resnick, S. (1984).
\newblock Limiting behaviour of sums and the term of maximum modulus.
\newblock {\em Proc. Lond. Math. Soc.}, {\em 49\/}(3), 385--422.

\bibitem[\protect\citeauthoryear{Pitman \& Yor}{Pitman \& Yor}{1996}]{PY1996}
Pitman, J. \& Yor, M. (1996).
\newblock Random discrete distributions derived from self-similar random sets.
\newblock {\em Electronic Journal of Probability}, {\em 1}, No. 4.

\bibitem[\protect\citeauthoryear{Pitman \& Yor}{Pitman \& Yor}{1997}]{PY1997}
Pitman, J. \& Yor, M. (1997).
\newblock The two-parameter {Poisson--Dirichlet} distribution derived from a
  stable subordinator.
\newblock {\em Annals of Probability}, {\em 25\/}(2), 855--900.

\bibitem[\protect\citeauthoryear{Resnick}{Resnick}{1986}]{resnick:1986b}
Resnick, S. (1986).
\newblock Point processes, regular variation and weak convergence.
\newblock {\em Adv. Appl. Probab.}, {\em 18}, 66--138.

\bibitem[\protect\citeauthoryear{Resnick}{Resnick}{1987}]{res87}
Resnick, S.~I. (1987).
\newblock {\em Extreme Values, Regular Variation, and Point Processes}.
\newblock Springer-Verlag.

\bibitem[\protect\citeauthoryear{Samorodnitsky \& Taqqu}{Samorodnitsky \&
  Taqqu}{1994}]{samtaqqu}
Samorodnitsky, G. \& Taqqu, M.~S. (1994).
\newblock {\em Stable non-{G}aussian random processes: stochastic models with
  infinite variance}.
\newblock Chapman \& Hall, London.

\bibitem[\protect\citeauthoryear{Shorrock}{Shorrock}{1972}]{shorrock1972}
Shorrock, R. (1972).
\newblock On record values and record times.
\newblock {\em J. Appl. Probab.}, {\em 3}, 316--326.

\bibitem[\protect\citeauthoryear{Smid \& Stam}{Smid \&
  Stam}{1975}]{smidstam1975}
Smid, B. \& Stam, A. (1975).
\newblock Convergence in distribution of quotients of order statistics.
\newblock {\em Stochastic Processes and their Applications}, {\em 3}, 287--292.

\bibitem[\protect\citeauthoryear{Teugels}{Teugels}{1981}]{teugels1981}
Teugels, J. (1981).
\newblock Limit theorems on order statistics.
\newblock {\em Ann. Probab.}, {\em 9}, 868--880.

\end{thebibliography}
\bibliographystyle{newapa}

\end{document}